\documentclass{article}
\usepackage{amsmath,amssymb}
\begin{document}
\begin{center}
{\Large BOUNDARY-VALUE PROBLEMS WITH NON-LOCAL INITIAL CONDITION FOR PARABOLIC EQUATIONS WITH PARAMETER}

\bigskip

{\bf J.M.Rassias$^1$ and E.T.Karimov$^2$}

\medskip

$^1${\it Pedagogical department, Mathematics and Informatics section, National and Capodostrian University of Athens.\,\,\,\,jrassias@primedu.uoa.gr}

$^2${\it Institute of Mathematics and Information Technologies, Uzbek Academy of Sciences.\,\,\,\,E-mail: erkinjon@gmail.com}

\end{center}

\begin{abstract}
In 2002, J.M.Rassias (Uniqueness of quasi-regular solutions for bi-parabolic elliptic bi-hyperbolic Tricomi problem, Complex Variables, 47 (8) (2002), 707-718) imposed and investigated the bi-parabolic elliptic bi-hyperbolic mixed type partial differential equation of second order. In the present paper some boundary-value problems with non-local initial condition for model and degenerate parabolic equations with parameter were considered. Also uniqueness theorems are proved and non-trivial solutions of certain non-local problems for forward-backward parabolic equation with parameter are investigated at specific values of this parameter by employing the classical "a-b-c" method. Classical references in this field of mixed type partial differential equations are given by: J.M.Rassias (Lecture Notes on Mixed Type Partial Differential Equations, World Scientific, 1990, pp.1-144) and M.M.Smirnov (Equations of Mixed Type, Translations of Mathematical Monographies, 51, American Mathematical Society, Providence, R.I., 1978 pp.1-232). Other investigations are achieved by G.C.Wen et al. (in period 1990-2007).
\end{abstract}

{\bf MSC 2000 classification:}


\section{Introduction}
Degenerate partial differential equations have numerous applications in Aerodynamics and Hydrodynamics. For example, problems for mixed subsonic and supersonic flows were considered by F.I.Frankl [1]. Reviews of interesting results on degenerated elliptic and hyperbolic equations up to 1965, one can find in the book by M.M.Smirnov [2]. Among other research results on this kind of equations were investigated by J.M.Rassias [3-10], G.C.Wen [11-15], A.Hasanov [16] and references therein. Also works by M.Gevrey [17], A.Friedman [18], Yu.Gorkov [19] are well-known on construction fundamental solutions for degenerated parabolic equations. In addition it was studied by N.N.Shopolov [20] a boundary-value problem with initial non-local condition for model parabolic equation in  [21-25]. However, in 2002,  J.M.Rassias (Uniqueness of quasi-regular solutions for bi-parabolic elliptic bi-hyperbolic Tricomi problem, Complex Variables, 47 (8) (2002), 707-718) imposed and investigated the bi-parabolic elliptic bi-hyperbolic mixed type partial differential equation of second order. In the present paper some boundary-value problems with non-local initial condition for model and degenerate parabolic equations with parameter were considered. Also uniqueness theorems are proved and non-trivial solutions of certain non-local problems for forward-backward parabolic equation with parameter are investigated at specific values of this parameter by employing the classical "a-b-c" method. Classical references in this field of mixed type partial differential equations are given by: J.M.Rassias (Lecture Notes on Mixed Type Partial Differential Equations, World Scientific, 1990, pp.1-144) and M.M.Smirnov (Equations of Mixed Type, Translations of Mathematical Monographies, 51, American Mathematical Society, Providence, R.I., 1978 pp.1-232). Other investigations are achieved by G.C.Wen et al. (in period 1990-2007).

\section{Non-local problems for degenerate parabolic equations with parameter.}
Let us consider a parabolic equation
$$
y^m u_{xx}  - x^n u_y  - \lambda x^n y^m u = 0, \eqno (1)
$$
with two lines of degeneration in the domain  $\Phi  = \left\{ {\left( {x,y} \right):\,0 < x < 1,\,\,0 < y < 1} \right\}$
, where $m,\,n > 0,\,\,\lambda  \in C $.

\textbf{The problem 1.} To find a regular solution of the equation satisfying boundary conditions
$$
u\left( {0,y} \right) = 0,\,\,\,\,\,u\left( {1,y} \right) = 0,\,\,\,\,0 \le y \le 1, \eqno (2)
$$
and non-local initial condition
$$
u\left( {x,0} \right) = \alpha u\left( {x,1} \right),\,\,0 \le x \le 1,\eqno (3)
$$
where $\alpha $ is non-zero real number.

The following statements are true:

\textbf{Theorem 1.} Let $\alpha  \in \left[ { - 1,0} \right) \cup \left( {0,1} \right]$
, ${\mathop{\rm Re}\nolimits} \lambda  \ge 0$. If there exists a solution of the problem 1, then it is unique.

\textbf{Corollary 1.} The problem 1 can have non-trivial solutions only when parameter $\lambda $
 lies outside of the sector $\Delta  = \left\{ {\lambda :\,\,{\mathop{\rm Re}\nolimits} \lambda  \ge 0} \right\}$
. These non-trivial solutions represented by

$$
u_{pk} \left( {x,y} \right) = C_{pk} \left( {\frac{2}{{n + 2}}} \right)^{\frac{1}{{n + 2}}} \mu _k^{\frac{1}{{2(n + 2)}}} x^{\frac{1}{2}} I_{\frac{1}{{n + 2}}} \left( {\frac{{2\sqrt {\mu _k } }}{{n + 2}}x^{\frac{{n + 2}}{2}} } \right)e^{\left( { - \ln |\alpha | - ip\pi } \right)y^{m + 1} }, \eqno (4)
$$
where $C_{pk} $ are constants, $p, k$ are real positive numbers. Eigenvalues defined as
$$
\lambda _{pk}  = \mu _k  + \left( {m + 1} \right)\ln |\alpha | + i\left( {m + 1} \right)p\pi.
$$
Here $\mu _k $ are roots of the equation
$$
I_{\frac{1}{{n + 2}}} \left( {\frac{{2\sqrt \mu  }}{{n + 2}}} \right) = 0,
$$
where $I_s \left(  \right)$ is the first kind modified Bessel function of $s$-th order.

We will omit the proof, because further we consider similar problem in three-dimensional domain in a full detail.

Let $\Omega $ be a simple-connected bounded domain in $R^3$ with boundaries $S_i$ ($i = \overline{1, 6}$). Here
$$
\begin{array}{l}
 S_1  = \left\{ {\left( {x,y,t} \right):\,t = 0,\,0 < x < 1,\,0 < y < 1} \right\},\,S_2  = \left\{ {\left( {x,y,t} \right):\,x = 1,\,0 < y < 1,\,0 < t < 1} \right\}, \\
 S_3  = \left\{ {\left( {x,y,t} \right):\,y = 0,\,0 < x < 1,\,0 < t < 1} \right\},\,\,S_4  = \left\{ {\left( {x,y,t} \right):\,x = 0,\,0 < y < 1,\,0 < t < 1} \right\}, \\
 S_5  = \left\{ {\left( {x,y,t} \right):\,y = 1,\,0 < x < 1,\,0 < t < 1} \right\},\,\,S_6  = \left\{ {\left( {x,y,t} \right):\,t = 1,\,0 < x < 1,\,0 < y < 1} \right\}. \\
 \end{array}
 $$

We consider the following degenerate parabolic equation
$$
x^n y^m u_t  = y^m u_{xx}  + x^n u_{yy}  - \lambda x^n y^m u \eqno (5)
$$
in the domain $\Omega $. Here $m > 0,\,n > 0,\,\lambda  = \lambda _1  + i\lambda _2 ,\,\,\lambda _1 ,\lambda _2  \in R$.

\textbf{The problem 2.} To find a function $u\left( {x,y,t} \right)$
 satisfying the following conditions:

i) $u\left( {x,y,t} \right) \in C\left( {\overline \Omega  } \right) \cap C_{x,y,t}^{2,2,1} \left( \Omega  \right)$;

ii) $u\left( {x,y,t} \right)$ satisfies the equation (5) in $\Omega $;

iii) $u\left( {x,y,t} \right)$ satisfies boundary conditions
$$
u\left( {x,y,t} \right)\left| {_{S_2  \cup S_3  \cup S_4  \cup S_5 }  = 0} \right.;\eqno (6)
$$

iv) and non-local initial condition
$$
u\left( {x,y,0} \right) = \alpha u\left( {x,y,1} \right).\eqno (7)
$$

Here $\alpha  = \alpha _1  + i\alpha _2 ,\,\,\alpha _1 ,\alpha _2 $ are real numbers, moreover $\alpha _1^2  + \alpha _2^2  \ne 0$.

\textbf{Theorem 2.} If $\alpha _1^2  + \alpha _2^2  < 1,\,\,\lambda _1  \ge 0$ and exists a solution of the problem 2, then it is unique.

\emph{Proof:}

Let us suppose that the problem 2 has two $u_1 ,\,u_2 $ solutions. Denoting $u = u_1  - u_2 $
 we claim that $u \equiv 0$ in $\Omega $.

First we multiply equation (5) to the function $\overline u \left( {x,y,t} \right)$, which is complex conjugate function of $u\left( {x,y,t} \right)$. Then integrate it along the domain $\Omega _\varepsilon$ with boundaries
$$
\begin{array}{l}
 S_{1\varepsilon }  = \left\{ {\left( {x,y,t} \right):\,t = \varepsilon ,\,\varepsilon  < x < 1 - \varepsilon ,\,\varepsilon  < y < 1 - \varepsilon } \right\},\\
 S_{2\varepsilon }  = \left\{ {\left( {x,y,t} \right):\,x = 1 - \varepsilon ,\,\varepsilon  < y < 1 - \varepsilon ,\,\varepsilon  < t < 1 - \varepsilon } \right\}, \\
 S_{3\varepsilon }  = \left\{ {\left( {x,y,t} \right):\,y = \varepsilon ,\,\varepsilon  < x < 1 - \varepsilon ,\,\varepsilon  < t < 1 - \varepsilon } \right\},\\
 S_{4\varepsilon }  = \left\{ {\left( {x,y,t} \right):\,x = \varepsilon ,\,\varepsilon  < y < 1 - \varepsilon ,\,\varepsilon  < t < 1 - \varepsilon } \right\}, \\
 S_{5\varepsilon }  = \left\{ {\left( {x,y,t} \right):\,y = 1 - \varepsilon ,\,\varepsilon  < x < 1 - \varepsilon ,\,\varepsilon  < t < 1 - \varepsilon } \right\},\\
 S_{6\varepsilon }  = \left\{ {\left( {x,y,t} \right):\,t = 1 - \varepsilon ,\,\varepsilon  < x < 1 - \varepsilon ,\,\varepsilon  < y < 1 - \varepsilon } \right\}. \\
 \end{array}
$$
Then taking real part of the obtained equality and considering
$$
{\mathop{\rm Re}\nolimits} \left( {y^m \overline u u_{xx} } \right) = {\mathop{\rm Re}\nolimits} \left( {y^m \overline u u_x } \right)_x  - y^m \left| {u_x } \right|^2 ,\,{\mathop{\rm Re}\nolimits} \left( {x^n \overline u u_{yy} } \right) = {\mathop{\rm Re}\nolimits} \left( {x^n \overline u u_y } \right)_y  - x^n \left| {u_y } \right|^2 ,
$$
$$
{\mathop{\rm Re}\nolimits} \left( {x^n y^m \overline u u_t } \right) = \left( {\frac{1}{2}x^n y^m \left| u \right|^2 } \right)_t ,
$$
after using Green's formula we pass to the limit at $\varepsilon  \to 0$. Then we get
$$
\begin{array}{l}
\int\limits_{\partial\Omega}\int{\mathop{\rm Re}\nolimits} \left[ {y^m \overline u u_x \cos \left( {\nu ,x} \right) + x^n \overline u u_y \cos \left( {\nu ,y} \right) - \frac{1}{2}x^n y^m \left| u \right|^2 \cos \left( {\nu ,t} \right)} \right]d\tau  \\
  = \int\int\limits_{\Omega}\int\left( {y^m \left| {u_x } \right|^2  + x^n \left| {u_y } \right|^2  + \lambda _1 x^n y^m \left| u \right|} \right)d\sigma  \\
 \end{array}
 $$
where $\nu $ is outer normal.
Taking into account
$
{\mathop{\rm Re}\nolimits} \left[ {\overline u u_x } \right] = {\mathop{\rm Re}\nolimits} \left[ {u\overline u _x } \right],\,\,{\mathop{\rm Re}\nolimits} \left[ {\overline u u_y } \right] = {\mathop{\rm Re}\nolimits} \left[ {u\overline u _y } \right]
$
we obtain
$$
{\mathop{\rm Re}\nolimits} \int {\int\limits_{S_1 } {\frac{1}{2}x^n y^m \left| u \right|^2 d\tau _1 } }  + \int {\int\limits_{S_2 } {y^m {\mathop{\rm Re}\nolimits} \left[ {u\overline u _x } \right]d\tau _2 } }  - \int {\int\limits_{S_3 } {x^n {\mathop{\rm Re}\nolimits} \left[ {u\overline u _y } \right]d\tau _3 } }  - \int {\int\limits_{S_4 } {y^m {\mathop{\rm Re}\nolimits} \left[ {u\overline u _x } \right]d\tau _4 } }  +
$$
$$
 + \int {\int\limits_{S_5 } {x^n {\mathop{\rm Re}\nolimits} \left[ {u\overline u _y } \right]d\tau _5 } }  - {\mathop{\rm Re}\nolimits} \int {\int\limits_{S_6 } {\frac{1}{2}x^n y^m \left| u \right|^2 d\tau _6  = } } \int {\int\limits_\Omega  {\int {\left( {y^m \left| {u_x } \right|^2  + x^n \left| {u_y } \right|^2  + \lambda _1 x^n y^m \left| u \right|} \right)d\sigma } } }.\eqno (8)
 $$
From (8) and by using conditions (6), (7), we find
$$
\frac{1}{2}\left[ {1 - \left( {\alpha _1^2  + \alpha _2^2 } \right)} \right]\int\limits_0^1 {\int\limits_0^1 {x^n y^m \left| {u\left( {x,y,1} \right)} \right|dxdy}  + \int {\int\limits_\Omega  {\int {\left( {y^m \left| {u_x } \right|^2  + x^n \left| {u_y } \right|^2  + \lambda _1 x^n y^m \left| u \right|} \right)d\sigma } } } }=0 .\eqno (9)
$$
Setting $\alpha _1^2  + \alpha _2^2  < 1,\,\,\lambda _1  \ge 0$, from (9) we have $u\left( {x,y,t} \right) \equiv 0$
in $\overline \Omega  $.

Theorem is proved.

We find below non-trivial solutions of the problem 2 at some values of parameter $\lambda $
 for which the uniqueness condition ${\mathop{\rm Re}\nolimits} \lambda  = \lambda _1  \ge 0$
 is not fulfilled.

We search the solution of Problem 2 as follows
$$
u\left( {x,y,t} \right) = X\left( x \right) \cdot Y\left( y \right) \cdot T\left( t \right).\eqno (10)
$$
After some evaluations we obtain the following eigenvalue problems:
$$
\left\{ \begin{array}{l}
 X''\left( x \right) + \mu _1 x^n X\left( x \right) = 0 \\
 X\left( 0 \right) = 0,\,\,\,X\left( 1 \right) = 0; \\
 \end{array} \right.\eqno (11)
 $$
$$
\left\{ \begin{array}{l}
 Y''\left( y \right) + \mu _2 y^m Y\left( y \right) = 0 \\
 Y\left( 0 \right) = 0,\,\,\,Y\left( 1 \right) = 0; \\
 \end{array} \right.\eqno (12)
 $$
$$
\left\{ \begin{array}{l}
 T'\left( t \right) + \left( {\lambda  + \mu } \right)T\left( t \right) = 0 \\
 T\left( 0 \right) = \alpha T\left( 1 \right). \\
 \end{array} \right.\eqno (13)
 $$
Here $\mu  = \mu _1  + \mu _2 $ is a Fourier constant.

Solving eigenvalue problems (11), (12) we find
$$
\mu _{1k}  = \left( {\frac{{n + 2}}{2}\widetilde{\mu _{1k} }} \right)^2 ,\,\,\,\,\mu _{2p}  = \left( {\frac{{m + 2}}{2}\widetilde{\mu _{2p} }} \right)^2 ,\eqno	(14)
$$
$$
X_k \left( x \right) = A_k \left( {\frac{2}{{n + 2}}} \right)^{\frac{1}{{n + 2}}} \mu _{1k}^{\frac{1}{{2\left( {n + 2} \right)}}} x^{\frac{1}{2}} J_{\frac{1}{{n + 2}}} \left( {\frac{{2\sqrt {\mu _{1k} } }}{{n + 2}}x^{\frac{{n + 2}}{2}} } \right),\eqno (15)
$$
$$
Y_p \left( y \right) = B_p \left( {\frac{2}{{m + 2}}} \right)^{\frac{1}{{m + 2}}} \mu _{2p}^{\frac{1}{{2\left( {m + 2} \right)}}} y^{\frac{1}{2}} J_{\frac{1}{{m + 2}}} \left( {\frac{{2\sqrt {\mu _{2p} } }}{{m + 2}}x^{\frac{{m + 2}}{2}} } \right),\eqno (16)
$$
where $k,p = 1,2,...$, $\widetilde{\mu _{1k} }$ and $\widetilde{\mu _{2p} }$ are roots of equations $J_{\frac{1}{{n + 2}}} \left( x \right) = 0$
 and $J_{\frac{1}{{m + 2}}} \left( y \right) = 0$, respectively.

The eigenvalue problem (13) has non-trivial solution only when
$\left\{ \begin{array}{l}
 \alpha _1  = e^{\lambda _1  + \mu _{kp} } \cos \lambda _2  \\
 \alpha _2  = e^{\lambda _1  + \mu _{kp} } \sin \lambda _2 . \\
 \end{array} \right.$

Here $\lambda  = \lambda _1  + i\lambda _2 ,\,\,\,\,\alpha  = \alpha _1  + i\alpha _2 ,\,\,\,\mu _{kp}  = \mu _{1k}  + \mu _{2p} $. After elementary calculations, we get
$$
\lambda _1  =  - \mu _{kp}  + \ln \sqrt {\alpha _1^2  + \alpha _2^2 } ,\,\,\,\,\lambda _2  = \arctan \frac{{\alpha _2 }}{{\alpha _1 }} + s\pi ,\,\,\,s \in Z^+\eqno (17)
$$
Corresponding eigenfunctions have the form
$$
T_{kp} \left( t \right) = C_{kp} e^{\left[ {\mu _{kp}  - \ln \sqrt {\alpha _1^2  + \alpha _2^2 }  - i\left( {\arctan \frac{{\alpha _2 }}{{\alpha _1 }} + s\pi } \right)} \right]t} .\eqno	(18)
$$

Considering (10), (15), (16) and (18) we can write non-trivial solutions of the problem 2 in the following form:
$$
u_{kp} \left( {x,y,t} \right) = D_{kp} \left( {\frac{2}{{n + 2}}} \right)^{\frac{1}{{n + 2}}} \left( {\frac{2}{{m + 2}}} \right)^{\frac{1}{{m + 2}}} \mu _{1k}^{\frac{1}{{2\left( {n + 2} \right)}}} \mu _{2p}^{\frac{1}{{2\left( {m + 2} \right)}}} \sqrt {xy} J_{\frac{1}{{n + 2}}} \left( {\frac{{2\sqrt {\mu _{1k} } }}{{n + 2}}x^{\frac{{n + 2}}{2}} } \right)$$

$$
 \times J_{\frac{1}{{m + 2}}} \left( {\frac{{2\sqrt {\mu _{2p} } }}{{m + 2}}y^{\frac{{m + 2}}{2}} } \right)
 e^{\left[ {\mu _{kp}  - \ln \sqrt {\alpha _1^2  + \alpha _2^2 }  - i\left( {\arctan \frac{{\alpha _2 }}{{\alpha _1 }} + s\pi } \right)} \right]t},
 $$ 									
where $D_{kp}  = A_k  \cdot B_p  \cdot C_{kp} $ are constants.

\textbf{Remark 1.} One can easily see that $\lambda _1  < 0$
 in (17), which contradicts to condition ${\mathop{\rm Re}\nolimits} \lambda  = \lambda _1  \ge 0$
 of the theorem 2.

\textbf{Remark 2.} The following problems can be studied by similar way.
Instead of condition (6) we put conditions as follows:

\begin{tabular}{|c|c|c|c|c|c|c|c|c|}
  \hline
  Problem's name & P$_3$ & P$_4$ & P$_5$ & P$_6$ & P$_7$ & P$_8$ & P$_9$ & P$_{10}$ \\
  \hline
  S$_2$ & $u_x$ & $u$ & $u$ & $u_x$ & $u$ & $u$ & $u_x$ & $u$ \\
  S$_3$ & $u_y$ & $u$ & $u_y$ & $u$ & $u_y$ & $u$ & $u$ & $u$ \\
  S$_4$ & $u$ & $u_x$ & $u$ & $u_x$ & $u$ & $u_x$ & $u$ & $u$ \\
  S$_5$ & $u$ & $u_y$ & $u_y$ & $u$ & $u$ & $u$ & $u$ & $u_y$ \\
  \hline
\end{tabular}

\section{Non-local problem for "forward-backward" parabolic equation with parameter.}
In the domain $D = D_1  \cup D_2  \cup I_0 ,$
$\,\,D_1  = \left\{ {\left( {x,y} \right): - 1 \le x \le 0,\,\,0 \le y \le 1} \right\},$

$I_0  = \left\{ {\left( {x,y} \right):\,x = 0,\,0 \le y \le 1\,} \right\}$
, $\,\,D_2  = \left\{ {\left( {x,y} \right):0 \le x \le 1,\,\,0 \le y \le 1} \right\}$
 let us  consider equation
$$Lu = \lambda u,\eqno (19)$$
where $\lambda  \in R,$ $Lu = u_{xx}  - sign\left( x \right)u_y $
.

\textbf{The problem 3.} To find a regular solution of the equation (19) from the class of functions $u\left( {x,y} \right) \in C\left( {\overline D } \right) \cap C^1 \left( {D \cup I_1  \cup I_2 } \right)$
, satisfying non-local conditions
$$k_1 u_x ( - 1,y) + k_2 u( - 1,y) = k_3 u_x (1,y),\,\,\,0 \le y \le 1,\eqno (20)$$
$$k_4 u_x (1,y) + k_5 u(1,y) = k_6 u_x ( - 1,y),\,\,\,0 \le y \le 1;\eqno (21)$$
$$u\left( {x,0} \right) = \alpha u\left( {x,1} \right),\,\,\,-1 \le x \le 1.\eqno (22)$$
Here $k_i \,\left( {i = \overline {1,6} } \right)$, $\alpha $
 is given non-zero constant, $I_1  = \left\{ {\left( {x,y} \right):\,x =  - 1,\,0 \le y \le 1\,} \right\}$,\\
  $I_2  = \left\{ {\left( {x,y} \right):\,x = 1,\,0 \le y \le 1\,} \right\}$.

Note, non-local conditions (20), (21) were used for the first time by N.I.Ionkin and E.I.Moiseev [26, 27].

\textbf{Theorem 3.} If
$$
|\alpha| =1,\,\,\lambda  > 0,\,k_3 k_5  = k_2 k_6 ,\,\,k_1 k_2  < 0,\,\,k_4 k_5  > 0\eqno (23)
$$
and exists a solution of the problem 3, then it is unique.

\emph{Proof:}

We multiply equation (19) to the function $u\left( {x,y} \right)$ and integrate along the domains $D_1 $ and $D_2 $.
Using Green's formula and condition (22), we get
$$
\int\limits_0^1 {u\left( { - 0,y} \right)u_x \left( { - 0,y} \right)dy}  = \int\limits_{ - 1}^0 {\frac{{\alpha ^2  - 1}}{2}u^2 \left( {x,1} \right)dx}  + \int\limits_0^1 {u\left( { - 1,y} \right)u_x \left( { - 1,y} \right)dy}  + \int\int\limits_{D_1}\left(u_x^2+\lambda u^2\right)dxdy,
$$
$$
\int\limits_0^1 {u\left( { + 0,y} \right)u_x \left( { + 0,y} \right)dy}  = \int\limits_0^1 {\frac{{\alpha ^2  - 1}}{2}u^2 \left( {x,1} \right)dx}  + \int\limits_0^1 {u\left( {1,y} \right)u_x \left( {1,y} \right)dy}  - \int\int\limits_{D_2}\left(u_x^2+\lambda u^2\right)dxdy.
$$
From conditions (20), (21), we find
$$
u\left( {1,y} \right)u_x \left( {1,y} \right) = \frac{{k_6 }}{{k_5 }}u_x \left( { - 1,y} \right)u_x \left( {1,y} \right) - \frac{{k_4 }}{{k_5 }}u_x^2 \left( {1,y} \right),
$$
$$
u\left( { - 1,y} \right)u_x \left( { - 1,y} \right) = \frac{{k_3 }}{{k_2 }}u_x \left( { - 1,y} \right)u_x \left( {1,y} \right) - \frac{{k_1 }}{{k_2 }}u_x^2 \left( { - 1,y} \right).
$$
Taking ainto account above identities we establish
$$
\begin{array}{l}
 \int\limits_{ - 1}^0 {\frac{{\alpha ^2  - 1}}{2}u^2 \left( {x,1} \right)dx}  + \int\limits_0^1 {\frac{{1-\alpha ^2}}{2}u^2 \left( {x,1} \right)dx}  + \int\limits_0^1 {\left[ {\frac{{k_4 }}{{k_5 }}u_x^2 \left( {1,y} \right) - \frac{{k_1 }}{{k_2 }}u_x^2 \left( { - 1,y} \right)} \right]dy + }  \\
  + \int\limits_0^1 {\left[ {\frac{{k_3 }}{{k_2 }} - \frac{{k_6 }}{{k_5 }}} \right]u\left( { - 1,y} \right)u_x \left( {1,y} \right)dy}  + \int\int\limits_{D_1}\left(u_x^2+\lambda u^2\right)dxdy + \int\int\limits_{D_2}\left(u_x^2+\lambda u^2\right)dxdy. \\
 \end{array}
$$

Considering condition (23), we get $u\left( {x,y} \right) \equiv 0$
 in $D$ an the proof of theorem is complete.

\textbf{Remark 3.} By similar method one can prove the uniqueness of solution of boundary-value problem with non-local initial condition for equation
$$
0 = \left\{ \begin{array}{l}
 y^m u_{xx}  + \left( { - x} \right)^n u_y  - \lambda \left( { - x} \right)^n y^m u = 0,\,\,\,\,x < 0 \\
 y^m u_{xx}  - x^n u_y  - \lambda x^n y^m u = 0,\,\,\,x > 0. \\
 \end{array} \right.
 $$

\textbf{Open question.} A question is still open, on the unique solvability of boundary value problems for the following equation:
$$
0 = \left\{ \begin{array}{l}
 y^{m_1 } \left( { - x} \right)^{n_2 } u_{xx}  + \left( { - x} \right)^{n_1 } y^{m_2 } u_y  - \lambda _1 u = 0,\,\,\,\,x < 0 \\
 y^{m_1 } x^{n_2 } u_{xx}  - x^{n_1 } y^{m_2 } u_y  - \lambda _2 u = 0,\,\,\,x > 0, \\
 \end{array} \right.
 $$
where $\lambda _1 ,\,\lambda _2 $ are given complex numbers and $m_i ,\,n_i  = const > 0\,\,\,\left( {i = 1,2} \right)$.

\begin{center}
\textbf{References}
\end{center}

1. Frankl F.I. On the problems of Chaplygin for mixed subsonic and supersonic flows. Izv.Akad. Nauk SSSR Ser. Mat., 1945, No 9, pp. 121-143.

2. Smirnov M.M. Degenerate elliptic and hyperbolic equations. Moscow: Nauka, 1966, 292 p.

3. J.M.Rassias, Uniqueness of Quasi-Regular Solutions for a Bi-Parabolic Elliptic Bi-Hyperbolic Tricomi Problem. Complex Variables, 2002, Vol.47, No 8, pp.707-718.

4. J.M.Rassias, Mixed Type Partial Differential Equations in R$^n$, Ph.D. Dissertation, university of California, Berkeley, USA, 1977, 1-135.

5. J.M.Rassias, Lecture Notes on Mixed Type Partial Differential Equations, World Scientific, 1990, pp.1-144.

6. J.M.Rassias, Mixed type partial differential equations with initial and boundary values in fluid mechanics. Int.J.Appl.Math.Stat. 13(2008), No J08, 77-107.

7. J.M.Rassias and A.Hasanov, Fundamental solutions of two degenerated elliptic equations and solutions of boundary value problems in infinite area. Int.J.Appl.Math.Stat. 8(2007), No M07. 87-95.

8. J.M.Rassias, Tricomi-Protter problem of nD mixed type equations. Int.J.Appl.Math.Stat. 8(2007), No M07. 76-86.

9.  J.M.Rassias, Existence of weak solutions for a parabolic elliptic-hyperbolic Tricomi problem. Tsukuba J.Math., 23(1999), no.1. 37-54.

10. J.M.Rassias, Uniqueness of quasi-regular solutions for a parabolic elliptic-hyperbolic Tricomi problem. Bull.Inst.Math.Acad.Sinica, 25(1997), no. 4. 277-287.

11. G.C.Wen, The Exterior Tricomi Problem for Generalized Mixed Equations with Parabolic Degeneracy. Acta Mathematics Sinica, English Series, Vol.22, No 5, pp. 1385-1398.

12. G.C.Wen and D.Chen, Discontinuous Riemann-Hilbert problems for quasilinear degenerate elliptic complex equations of first order. Complex Var. Theory Appl. 50(2005), no. 7-11. 707-718.

13. G.C.Wen, The mixed boundary-value problem for second order elliptic equations with degenerate curve on the sides of an angle. Math. Nachr. 279(2006), no.13-14. 1602-1613.

14. S.Huang, Y.Y.Qiao and G.C.Wen, Real and complex Clifford analysis. Advances in Complex Analysis and its Applications, 5. Springer, New York, 2006. x+251.

15. G.C.Wen and H.G.W.Begehr, Boundary value problems for elliptic equations and systems. Pitman Monographs and Surveys in Pure and Applied Mathematics, 46. Longman Scientific and Tech., Harlow; John Wiley and Sons, Inc.,N.Y., 1990. xii+411.

16. A.Hasanov, Fundamental solutions of generalized bi-axially  symmetric Helmholtz equation.
Complex Variables, 2007, Vol. 52, No. 8, pp. 673-683.

17. M.Gevrey, Sur les equations aux derivees partielles du type parabolique. J.Math.Appl., 1914, Ch.4, pp.105-137.

18. A.Friedman, Fundamental solutions for degenerate parabolic equations. Acta Mathematica 133. Imprime le 18 Fevrier 1975, pp. 171-217.

19. Yu.P.Gorkov, Construction of a fundamental solution of parabolic equation with degeneration. Calcul. methods and programming, 2005, V.6, pp.66-70.

20. N.N.Shopolov, Mixed problem with non-local initial condition for a heat conduction equation. Reports of Bulgarian Academy of Sciences, 1981, V.34, No 7, pp.935-936.

21. A.A.Kerefov, On a spectra of the Gevrey problem for loaded mixed parabolic equation. Proceeding of the conference:Classical and non-classical boundary-value problems for partial differential equations, special functions, integral equations and their applications, 25-29 April, Russia, pp.76-77.

22. K.B.Sabitov, To the theory of mixed parabolic-hyperbolic type equations with spectral parameter. Differencial'nye Uravnenija, 1989, Vol.25, No 1, pp. 117-126.

23. A.S.Berdyshev and E.T.Karimov, Some non-local problems for the parabolic-hyperbolic type equation with non-characteristic line of changing type. CEJM, 2006, Vol. 4, No 2, pp. 183-193.

24. E.T.Karimov, About the Tricomi problem for the mixed parabolic-hyperbolic type equation with complex spectral parameter. Complex Variables, 2005, Vol.56, No 6, pp.433-440.

25. E.T.Karimov, Some non-local problems for the parabolic-hyperbolic type equation with complex spectral parameter. Mathematische Nachrichten, Vol.281, Iss. 7, 2008, pp.959-970.

26. N.I.Ionkin, The stability of a problem in the theory of heat condition with non-classical boundary conditions. (Russian). Differencial'nye Uravnenija 15 (1979), No 7, 1279-1283, 1343.

27. N.I.Ionkin and E.I.Moiseev, A problem for a heat equation with two-point boundary conditions. (Russian). Differencial'nye Uravnenija 15 (1979), No 7, 1284-1295, 1343.

\medskip

\textbf{Name: }JOHN MICHAEL RASSIAS

\textbf{Affiliation:} Professor

\textbf{Institute:} National and Capodistrian University of Athens

\textbf{Postal address:} 4, Agamemnonos Str., Aghia Paraskevi, Athens 15342, Greece

\medskip

\textbf{Name: }KARIMOV ERKINJON TULKINOVICH

\textbf{Affiliation:} Junior Scientific Researcher

\textbf{Institute:} Institute of Mathematics and information
technologies, Uzbek Academy of Sciences

\textbf{Postal address:} 29, Durmon yuli str., Tashkent-100125, Uzbekistan.

\end{document}